\documentclass[a4paper,10pt]{article}
\usepackage{amsmath}
\usepackage{amssymb,amsfonts,amsthm}
\usepackage[applemac]{inputenc}
\usepackage[english]{babel}
\usepackage[pdftex]{graphicx}
\usepackage{caption}
\usepackage{subcaption}
\usepackage{anysize}
\usepackage{enumerate}
\usepackage{tikz}
\usepackage[all]{xy}
\usepackage{mathrsfs}
\usepackage{verbatim}
\usepackage{cite}
\usepackage{functan} 
\definecolor{green}{rgb}{0,0.8,0.5}

\newcommand{\abs}[1]{\left\vert#1\right\vert}
\newcommand{\overbar}[1]{\mkern 1.5mu\overline{\mkern-1.5mu#1\mkern-1.5mu}\mkern 1.5mu}
\graphicspath{{/pictures}}

\usepackage[pdftex]{hyperref}
\hypersetup{colorlinks=true,linkcolor=blue,citecolor=red}
\usepackage[width=0.8\textwidth]{caption}

\renewenvironment{abstract}{\small\quotation\noindent
 {\bfseries \abstractname .}}{\endquotation \par}

\newenvironment{trivlistparm}[2]{\trivlist
\item[{#1 #2}]}{\endtrivlist}

\newenvironment{prooftext}[1]{\trivlistparm{\bfseries}{#1}}{\Qed\endtrivlistparm}
\newenvironment{demo}{\trivlistparm{\bfseries}{Proof.}}{\Qed\endtrivlistparm}

\catcode`\@=11

\def\resetthefootnote{\renewcommand{\thefootnote}{\@arabic\c@footnote} }
\def\@principiremex#1{\trivlist
 \item[\hskip \labelsep{\bfseries #1\ \thethm.}]\ignorespaces}
\def\opar@principiremex#1[#2]{\trivlist
 \item[\hskip \labelsep{\bfseries #1\ \thethm\ (#2).}]\ignorespaces}

\newcommand{\newTHEOremrom}[2]{\newenvironment{#1}{\refstepcounter{thm}\@ifnextchar[{\opar@principiremex{#2}}
{\@principiremex{#2}}}{\qedB\endtrivlist}} \catcode`\@=12
\DeclareMathSymbol{\square}{\mathord}{AMSa}{"03}
\newcommand{\qedB}{\nopagebreak\hspace*{\fill}$\square$\par}
\newcommand{\Qed}{\nopagebreak\hspace*{\fill}{\vrule width6pt height6pt depth0pt}\par}

\newTHEOremrom{defi} {Definition}
\newTHEOremrom {rem} {Remark}
\newTHEOremrom {ex} {Example}

\renewcommand{\geq}{\geqslant}
\renewcommand{\epsilon}{\varepsilon}
\renewcommand{\leq}{\leqslant}
\newcommand{\R}{\mathbb{R}}
\newcommand{\Q}{\mathbb{Q}}
\newcommand{\Z}{\mathbb{Z}}
\newcommand{\N}{\mathbb{N}}

\newcommand{\C}{\mathbb{C}}

\newtheorem{thm}{Theorem}[section]

\newtheorem{cor}[thm]{Corollary}
\newtheorem{lema}[thm]{Lemma}

\newtheorem{claim}{Claim}

\title{\textbf{A proof of Bertrand's theorem using the theory of isochronous potentials}
\footnotetext{All the authors are partially supported by the MINECO/FEDER grant MTM2017-82348-C2-1-P. D. Rojas is also partially supported by the MINECO/FEDER grant MTM2017-86795-C3-1-P.}
\footnotetext{{\it Email addresses:} \texttt{rojas@ugr.es}  (D.~Rojas, corresponding author), \texttt{rortega@ugr.es} (R.~Ortega).}
}

\author{Rafael Ortega and David Rojas \\[10pt]
{\small \textsl{Departamento de Matem\'atica Aplicada,}}\\
\vspace{-2pt}
{\small \textsl{Universidad de Granada, 18071 Granada, Spain}}}

\date{}

\begin{document}

\maketitle

\begin{abstract}
We give an alternative proof for the celebrated Bertrand's theorem as a corollary of the isochronicity of a certain family of centers.
\end{abstract}

\section{Introduction and main results}

Given a field of forces in the Euclidean space which is central and attractive there always exist circular periodic solutions. In $1873$ Bertrand~\cite{Bertrand} proved the following result: among all central fields of forces in the Euclidean space there are only two exceptional cases (the harmonic oscillator and the Newtonian potential) in which all solutions close to the circular motions are also periodic. From a philosophical point of view this result has strong implications on the Gravitation Law. Besides the original proof, nowadays there are several methods of proof. See~\cite{Albouy,Fejoz,Rusos} for more information.

The goal of the present paper is to show the connection of Bertrand's theorem with the theory of planar isochronous potential centers. Given a $C^2$-potential $V=V(x)$ defined on an open interval $I\subset\R$ and an equilibrium $x_*\in I$, $V'(x_*)=0$, we say that the potential $V$ is $T$-isochronous around $x_*$ if there exists a neighborhood $\mathcal U\subset\R^2$ of $(x_*,0)$ such that every solution of the system
\begin{equation}\label{sys}	
\dot x=y,\ \dot y=-V'(x)
\end{equation}
passing through $\mathcal U\setminus\{(x_*,0)\}$ is periodic with $T$ as minimal period. This is a local notion depending only of the behavior of $V$ around $x_*$. Isochronous potentials have been considered by many authors and we refer to~\cite{CMV,Urabe61} for more information. Next we introduce a notion of common isochronicity for families of potentials. Assume that $\lambda>0$ is a parameter and $V_{\lambda}=V_{\lambda}(x)$ is a $C^2$-potential defined on $I$. We say that the family $\{V_{\lambda}\}_{\lambda>0}$ is $T$-isochronous if each $V_{\lambda}$ has at most one equilibrium and the potential $V_{\lambda}$ is $T$-isochronous around it. To avoid trivial situations we also assume that $V_{\lambda}$ has an equilibrium for some $\lambda>0$. This definition can be illustrated by the simple families
\[
V_{\lambda}(x)=\frac{x^2}{2}+\lambda x \text{ and }\tilde V_{\lambda}(x)=\frac{\lambda x^2}{2},
\]
defined on the whole real line $I=\R$. The first family is $2\pi$-isochronous while the second is not $T$-isochronous for any $T>0$. Note that $\tilde V_{\lambda}$ is isochronous but the period depends upon $\lambda$.

We will be interested in families of potentials of the type
\begin{equation}\label{bertie}
V_{\lambda}(x)=\frac{x^2}{2}+\lambda\Phi(x),\ x\in(0,\infty)
\end{equation}
where $\Phi\in C^2(0,\infty)$ is a given function. In the next result we show that there are only two isochronous families of this type. The first family corresponds to a harmonic oscillator with constant forcing and the second appears frequently in the literature associated to different names (see~\cite{Er1,Er2,Pinney}).

\begin{thm}\label{fam}
There are only two families of the type~\eqref{bertie} which are $T$-isochronous for some $T>0$. They correspond to
\[
\Phi(x)=Kx \text{ and }\Phi(x)=-\frac{K}{4}x^{-4}
\]
for some $K<0$. In the first case $T=2\pi$ and in the second $T=\pi$.
\end{thm}

This result is somehow related to Theorems~$2$ and~$8$ in~\cite{Rusos}. There the notion of isochronous family is replaced by the stronger concept of rationally closing $\Phi'$.

We claim that Bertrand's theorem follows as a corollary of the previous result. The assumption on the smoothness of the central force will be just $C^1$. In previous proofs it has been assumed that the force was analytic or $C^{\infty}$. 

\begin{thm}[Bertrand's theorem]\label{bertrand}
Consider the central force problem
\begin{equation}\label{central}
\ddot x = -\phi(r^2)x,\  \ddot y=-\phi(r^2)y
\end{equation}
with $r^2=x^2+y^2$ and $\phi\in C^1(0,+\infty)$. Assume that the force is attractive somewhere and that all solutions near circular motions are periodic. Then there is a constant $K<0$ such that either the central force is $\phi(\eta)=K$ or $\phi(\eta)=K\eta^{-\frac{3}{2}}$.
\end{thm}

We refer to Section~\ref{sec:demo} for a more precise statement.

The rest of the paper is divided in three sections. In Section~\ref{sec:2} we review the elegant theory of isochronous potentials developed by Urabe in~\cite{Urabe61} and~\cite{Urabe62}. The original formulation was only valid in some small neighborhood of the equilibrium and we will work in a more global setting. Urabe's theory will allow us to derive a useful consequence on the behavior of the potential at the end points of the region of isochronicity. These ideas are applied in Section~\ref{sec:3} to prove Theorem~\ref{fam}. In Section~\ref{sec:demo} we prove Bertrand's Theorem as a consequence of Theorem~\ref{fam}. The proof is more or less direct when $\phi$ is real analytic but there are some subtleties when $\phi$ is only $C^1$.

\section{Remarks on Urabe's theorem}\label{sec:2}

Given an interval $I=(\alpha,\beta)$ with $-\infty\leq \alpha<0<\beta\leq +\infty$, the Urabe class $\mathcal U(I)$ is defined as the set of functions $u:I\rightarrow\R$ satisfying $u(0)=0$, $u\in C(I)$ and $v\in C^1(I)$ where $v(x)=xu(x)$. As an example we consider the function $u(x)=x^{\frac{1}{3}}$, it belongs to $\mathcal U(\R)$. It is not hard to prove that a function $u$ belongs to $\mathcal U(I)$ if and only if $u\in C(I)\cap C^1(I\setminus\{0\})$, $u(0)=0$ and $\lim_{x\rightarrow 0}x u'(x)=0$.

\begin{thm}\label{thm:urabe-global}
Assume that $V\in C^2(I)$ satisfies $V(0)=V'(0)=0$, $xV'(x)>0$ if $x\neq 0$ and there exists $\overbar{V}\in(0,+\infty]$ such that
\[
\lim_{x\rightarrow\alpha^+}V(x)=\lim_{x\rightarrow\beta^-} V(x)=\overbar{V}.
\]
In addition, assume that there exists $T>0$ such that for each $x_0\in I$, $x_0\neq 0$, the Cauchy problem
\[
\ddot x + V'(x)=0,\ x(0)=x_0,\ \dot x(0)=0
\]
has a solution with minimal period $T$. 

Then there exists an odd function $S\in\mathcal U(J)$ with $J=(-(2\overbar{V})^{\frac{1}{2}},+(2\overbar{V})^{\frac{1}{2}})$ and $\abs{S(X)}<1$ if $X\in J$ such that the solution $X(x)$ of
\begin{equation}\label{ivp}
\frac{dX}{dx}=\frac{2\pi}{T}\frac{1}{1+S(X)},\ X(0)=0
\end{equation}
is defined on the interval $I$ where it satisfies
\begin{equation}\label{sq}
V(x)=\frac{1}{2}X(x)^2.
\end{equation}
\end{thm}

\begin{rem}\label{rem:urabe}
\begin{enumerate}
\item The solution of~\eqref{ivp} is unique. In general the solution of the initial value problem $\frac{dX}{dx}=G(X)$ is unique as soon as $G$ is continuous and $G\neq 0$ everywhere. See~\cite{Corduneanu}.
\item The converse of the above Theorem also holds. Given a function $S(X)$ in the above conditions and a number $T>0$, the function $V(x)$ can be defined by the identities~\eqref{ivp} and~\eqref{sq}. It can be checked that this function satisfies all the conditions imposed in the direct Theorem. Note that $I=(\alpha,\beta)$ now is the maximal interval for the solution of~\eqref{ivp} and $X(x)\rightarrow -(2\overbar{V})^{\frac{1}{2}}$ as $x\rightarrow \alpha^+$ and $X(x)\rightarrow +(2\overbar{V})^{\frac{1}{2}}$ as $x\rightarrow\beta^-$.
\item An interesting consequence of this Theorem is the non-existence of potential isochronous centers whose period annulus is a vertical strip of the type $(\alpha,\beta)\times\R$ with $-\infty<\alpha<\beta<+\infty$. Assume by contradiction that every non-constant orbit with initial condition in $(\alpha,\beta)\times\R$ has minimal period $T$. The orbits approaching $(\alpha,0)$ and $(\beta,0)$ would cross the $y$-axis at points $(0,y)$ with $\abs{y}\rightarrow+\infty$. This would imply $\overbar{V}=+\infty$. We claim that the function $1+S$ should be integrable in the Lebesgue sense on the whole real line, that is $1+S\in L^1(\R)$. Since it is a positive function it is enough to prove that the integral is finite and this is a consequence of~\eqref{ivp}. Indeed,
\[
\int_0^{+\infty}(1+S(X))dX=\frac{2\pi}{T}\beta,\  \int_{-\infty}^0(1+S(X))dX=-\frac{2\pi}{T}\alpha.
\]
Since $S$ is odd,
\[
\int_0^{+\infty}(1-S(X))dX=-\frac{2\pi}{T}\alpha.
\]
From $\abs{S(X)}<1$ we deduce that $1-S(X)$ is positive and $1-S\in L^1(0,\infty)$. Summing up, both functions $1+S$ and $1-S$ belong to $L^1(0,\infty)$ and this is absurd because the sum function is the constant $2$, a non-integrable function.
\end{enumerate}\vspace*{-.55cm}
\end{rem}

The proof of Theorem~\ref{thm:urabe-global} requires two previous lemmas. The proof of the first one is immediate whereas the second is a consequence of~\cite{Urabe62}.

\begin{lema}\label{lema1}
Assume that $V\in C^2(I)$ satisfies $V(0)=V'(0)=0$, $xV'(x)>0$ if $x\neq 0$, $V''(0)\neq 0$ and 
\[
\lim_{x\rightarrow\alpha^+}V(x)=\lim_{x\rightarrow\beta^-}V(x)=\overbar{V}
\]
for some $\overbar{V}\in(0,+\infty]$. Then the function $X(x)\!:=\text{sign}(x)\sqrt{2V(x)}$ is a $C^1$ diffeomorphism from $I$ to $J$. Moreover, $X'(0)=+\sqrt{V''(0)}$.
\end{lema}

\begin{lema}\label{lema2}
Let $R>0$ and assume that $T\in C(-R,R)$ is even and
\[
\int_0^{2\pi} T(r\cos\theta)d\theta =0 \text{ for each }r\in(0,R).
\]
Then $T$ is zero everywhere.
\end{lema}

\begin{prooftext}{Proof of Theorem~\ref{thm:urabe-global}.}
By assumption the origin is an isochronous center with minimal period $T$ and this implies
\[
V''(0)=\left(\frac{2\pi}{T}\right)^2.
\]
Let us define $X(x)$ by Lemma~\ref{lema1}. We consider the map $\mathscr A\rightarrow\mathscr B$, $(x,y)\mapsto (X,Y)$, $X=X(x)$, $Y=y$ between the domains
\[
\mathscr A=\{(x,y)\in I\times\R : \frac{1}{2}y^2+V(x)<\overbar{V}\} \text{ and }
\mathscr B=\{(X,Y)\in J\times\R : Y^2+X^2<2\overbar{V}\}.
\]
This map is a $C^1$-diffeomorphism transporting those orbits of $\dot x=y$, $\dot y=-V(x)$ crossing the $x$-axis to concentric circles. The case $\overbar{V}<+\infty$ is illustrated in Figure~\ref{figura}. The system becomes in the new variables
\[
\begin{cases}
\dot X &= \frac{H(X)}{X}Y,\\
\dot Y &= -H(X),
\end{cases}
\]
where $H=V'\circ x$ and $x=x(X)$ is the inverse function of $X=X(x)$.
The associated vector field is continuous on $\mathscr B$ and there is uniqueness for the initial value problem because the system is equivalent to the system $\dot x=y$, $\dot y=-V(x)$. In addition, $X^2+Y^2$ is a first integral and the punctured disk $\mathscr B\setminus\{(0,0)\}$ is invariant. Using polar coordinates $X=r\cos\theta$, $Y=r\sin\theta$ the system is written as
\[
\begin{cases}
\dot r &= 0,\\
\dot \theta &= -\omega(r\cos\theta),
\end{cases}
\]
where $\omega(X)=\frac{H(X)}{X}$. We observe that $\omega(0)=\sqrt{V''(0)}=\frac{2\pi}{T}$. At this point the reader who is familiar with the theory of integrable Hamiltonian systems will realize that these coordinates designed by Urabe are reminiscent of the classical action-angle variables. However they are not the same, note that the above system in $(r,\theta)$ is not of Hamiltonian type.

\begin{figure}
\centering
\includegraphics[scale=1]{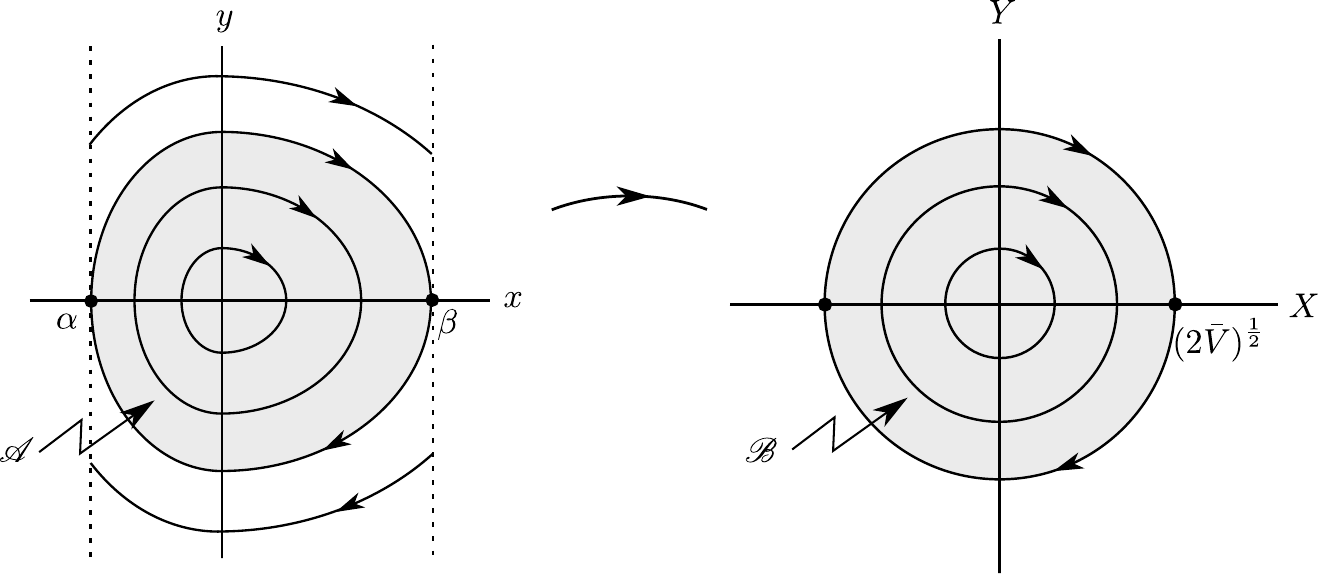}
\caption{\label{figura}Urabe's map.}
\end{figure}

The function $\omega(X)$ is positive on $J$ and all orbits on $\mathscr B\setminus\{(0,0)\}$ have minimal period $T$, this implies the identity
\begin{equation}\label{eqT}
T=\int_0^{2\pi}\frac{d\theta}{\omega(r\cos\theta)}
\end{equation}
for each $r\in(0,+(2\overbar{V})^{\frac{1}{2}})$. On the other hand, the function $u(X)=\frac{1}{\omega(X)}-\frac{T}{2\pi}$ belongs to $\mathcal U(J)$. It is clear that $u(0)=0$ and $u\in C(J)\cap C^1(J\setminus\{0\})$ and it remains to prove that $\lim_{X\rightarrow 0}Xu'(X)=0$. For $X\neq 0$,
\[
Xu'(X)=-X\frac{\omega'(X)}{\omega(X)^2}=\frac{1}{\omega(X)^2}\left(\frac{H(X)}{X}-H'(X)\right).
\]
By L'H\^opital rule, $\lim_{X\rightarrow 0}\frac{H(X)}{X}=\lim_{X\rightarrow 0} H'(X)=\frac{2\pi}{T}$, and then $Xu'(X)\rightarrow 0$ as $X\rightarrow 0$.

We can decompose $u(X)$ in the odd and even part so that 
\[
\frac{1}{\omega(X)}=\frac{T}{2\pi}\bigl(1+S(X)+T(X)\bigr)
\]
where $S(X)=\frac{\pi}{T}(u(X)-u(-X))$ and $T(X)=\frac{\pi}{T}(u(X)+u(-X))$. From these formulas we deduce that the functions $S(X)$ and $T(X)$ belong to $\mathcal U(J)$. Let us now prove that $T(X)$ is identically zero. The identity~\eqref{eqT} can be reformulated as
\[
\sigma(r)+\tau(r)=0,\ r\in J
\]
where
\[
\sigma(r)\!:=\int_0^{2\pi} S(r\cos\theta)d\theta \text{ and }
\tau(r)\!:=\int_0^{2\pi} T(r\cos\theta)d\theta.
\]
Since $\sigma$ is an odd function and $\tau$ is even we deduce that $\tau(r)=0$ for every $r\in J$. Then Lemma~\ref{lema2} implies that $T\equiv 0$. In consequence
\[
\omega(X)=\frac{2\pi}{T}\frac{1}{1+S(X)}.
\]
From the definition of $X(x)$ we deduce that~\eqref{ivp} and~\eqref{sq} hold. It remains to prove that $\abs{S(X)}<1$. We already know that $1+S(X)=\frac{2\pi}{T\omega(X)}>0$. Since $S$ is odd also $1-S(X)>0$. and the conclusion follows.
\end{prooftext}

We finish this Section with a result on the behavior of the potential at the end points of the interval of isochronicity.

\begin{cor}\label{uno}
Let $V$ be a potential in the conditions of Theorem~\ref{thm:urabe-global}. Then
\[
\limsup_{x\rightarrow\alpha^+}\frac{\sqrt{2V(x)}}{\abs{V'(x)}}+\liminf_{x\rightarrow\beta^-}\frac{\sqrt{2V(x)}}{V'(x)}=
\liminf_{x\rightarrow\alpha^+}\frac{\sqrt{2V(x)}}{\abs{V'(x)}}+\limsup_{x\rightarrow\beta^-}\frac{\sqrt{2V(x)}}{V'(x)}=\frac{T}{\pi}.
\]
\end{cor}

\begin{demo}
The function $S(X)$ is determined from the formula
\begin{equation}\label{ide3}
S(X(x))=-1+\frac{2\pi}{T}\frac{\sqrt{2V(x)}}{V'(x)}\mathrm{sign}(x).
\end{equation}
This is a consequence of~\eqref{ivp} and~\eqref{sq}. Also, from~\eqref{sq} we deduce that $X(x)\rightarrow -\sqrt{2\overbar{V}}$ as $x\rightarrow\alpha^+$ and $X(x)\rightarrow+\sqrt{2\overbar{V}}$ as $x\rightarrow\beta^-$. Since $S$ is odd, 
\[
\limsup_{X\rightarrow+\sqrt{2\overbar{V}}}S(X)=-\liminf_{X\rightarrow-\sqrt{2\overbar{V}}}S(X).
\] 
Obviously the $\sup$ and $\inf$ can be interchanged. The conclusion follows by letting $x$ to go to the end points $\alpha$ and $\beta$ in the identity~\eqref{ide3}.
\end{demo}

\section{Proof of Theorem~\ref{fam}}\label{sec:3}

Assume that $V=V(x)$ is a $C^2$ function defined on a neighborhood of $x_*$ with $V'(x_*)=0$. Then $x=x_*$ is an equilibrium of $\ddot x + V'(x)=0$ and it is well known that when $x_*$ is an isochronous center then $V''(x_*)$ is positive and the minimal period is given by the formula
\begin{equation}\label{T2}
T=\frac{2\pi}{\sqrt{V''(x_*)}}.
\end{equation}
This is a reformulation of a fact already mentioned in the proof of Theorem~\ref{thm:urabe-global}. Let us assume that $V_{\lambda}$ is a $T$-isochronous family of the type~\eqref{bertie}. Define the set
\[
\Lambda\!:=\{\lambda\in(0,\infty): V_{\lambda}\text{ has a critical point}\}.
\]
By assumption we know that $\Lambda$ is non-empty. From now on we employ the notation $\varphi=\Phi'$ so that $x(\lambda)$ is the unique solution of
\begin{equation}\label{imp}
x+\lambda\varphi(x)=0,\ x\in(0,\infty)
\end{equation}
when $\lambda\in\Lambda$. We claim that $\Lambda$ is an interval. Given $\lambda_1,\lambda_2\in\Lambda$, $\lambda_1<\lambda_2$, we take any $\lambda\in(\lambda_1,\lambda_2)$ and consider the continuous function $f(x)\!:=\frac{1}{\lambda}x+\varphi(x)$. Then $f(x(\lambda_1))<0<f(x(\lambda_2))$ and therefore $f$ has a zero between $x(\lambda_1)$ and $x(\lambda_2)$. In consequence $\lambda\in\Lambda$ and $\Lambda$ is an interval. We know by assumption that $V_{\lambda}$ is $T$-isochronous around $x(\lambda)$ and therefore $V_{\lambda}''(x(\lambda))=1+\lambda\varphi'(x(\lambda))>0$. The function $x=x(\lambda)$ is defined by the equation~\eqref{imp} and the implicit function theorem can be applied. We deduce that $\Lambda$ is open and $\lambda\in\Lambda\mapsto x(\lambda)\in(0,\infty)$ is a $C^1$ function with positive derivative everywhere. In consequence the image 
\[
\mathcal J\!:=\{ x(\lambda):\lambda\in\Lambda\}
\]
is an open interval and the inverse function $x\in\mathcal J\mapsto \lambda(x)\in\Lambda$ is an increasing diffeomorphism. The identity~\eqref{T2} together with the definition of $x(\lambda)$ yields to the differential equation
\begin{equation}\label{formula}
\varphi'(x)=\frac{a}{x}\varphi(x),\ x\in\mathcal J
\end{equation}
where $a\!:=1-\left(\frac{2\pi}{T}\right)^2$. We also notice that $\varphi(x)$ must be negative on $\mathcal J$. Solving the linear equation~\eqref{formula} we obtain $\varphi(x)=Kx^{a}$, $x\in\mathcal J$, for some negative $K$. The function $\lambda(x)$ is now easily computed,
\[
\lambda(x)=-\frac{x}{\varphi(x)}=-\frac{1}{K}x^{1-a},\ x\in\mathcal J.
\]
We claim that $\mathcal J=(0,\infty)$. Otherwise let $x_*\in\partial\mathcal J\cap(0,\infty)$. Then a sequence $x_n\in\mathcal J$ with $x_n\rightarrow x_*$ produces a sequence $\lambda_n=-\frac{1}{K}x_n^{1-a}\rightarrow -\frac{1}{K}x_*^{1-a}=\lambda_*\in(0,\infty)$ and $\lambda_*\in\Lambda$. Thus $x_*\in\mathcal J$ arriving to contradiction.

From the previous discussion we can write $V_{\lambda}(x)=\frac{x^2}{2}+\lambda K \frac{x^{a+1}}{a+1}$ with $a<1$ and $K<0$. We know that $V_{\lambda}$ is $T$-isochronous around $x(\lambda)$ and therefore $V(x)=\frac{x^2}{2}-\frac{x^{a+1}}{a+1}$, $a\neq -1$, or $V(x)=\frac{x^2}{2}-\ln x$, $a=-1$, must be $T$-isochronous around $x=1$. This is a consequence of the homogeneity of the function $x^a$ and a rescaling argument. It will be useful to simplify the computations. We intend to apply Corollary~\ref{uno} to this potential but then we must translate the equilibrium to the origin. From now on we consider the potential
\[
V(x)=\frac{(x+1)^2}{2}-\frac{(x+1)^{a+1}}{a+1}-\frac{1}{2}+\frac{1}{a+1},\ a\neq -1
\]
and
\[
V(x)=\frac{(x+1)^2}{2}-\ln(x+1)-\frac{1}{2},\ a= -1
\]
and assume that $V$ is $T$-isochronous around $x=0$. To adjust to the setting of Corollary~\ref{uno} we define $\alpha=-1$ and 
\[
\beta=\begin{cases}
-1+\left(\frac{2}{1+a}\right)^{\frac{1}{1-a}} & \text{ if }\abs{a}<1,\\
+\infty & \text{ if }a\leq -1.
\end{cases}
\]
We know that the solution of $\ddot x+V'(x)=0$, $x(0)=x_0$, $\dot x(0)=0$ has minimal period $T$ if $x_0\neq 0$ is sufficiently small. The function $x^a$ is analytic on the interval $(\alpha,\beta)$ and all the solutions of the above initial value problem with $x_0\in(\alpha,\beta)$ are periodic. Let $\tau=\tau(x_0)$ be the corresponding minimal period. This function is constant ($\tau\equiv T$) in a neighborhood of $x_0=0$ and analytic on $(\alpha,0)\cup(0,\beta)$. Therefore $\tau\equiv T$ everywhere and the assumptions of Theorem~\ref{thm:urabe-global} hold. In consequence the conclusion of Corollary~\ref{uno} holds. The isochronous potentials postulated by the Theorem correspond to $a=0$ and $a=-3$. We must exclude the remaining cases and to this end we distinguish cases. Assume first that $a\in(0,1)$. Then $\lim_{x\rightarrow \alpha^+} V'(x)=0$ and $\lim_{x\rightarrow\alpha^+}\frac{\sqrt{2V(x)}}{\abs{V'(x)}}=+\infty$. This is absurd because this last limit is below $\frac{T}{\pi}$. Next we assume $a\in(-1,0)$. In this case $\lim_{x\rightarrow\alpha^+}V'(x)=-\infty$ and $\lim_{x\rightarrow\alpha^+}\frac{\sqrt{2V(x)}}{\abs{V'(x)}}=0$. Then $\lim_{x\rightarrow\beta^-}\frac{\sqrt{2V(x)}}{V'(x)}=\frac{T}{\pi}$. Since $V(\beta)=-\frac{1}{2}+\frac{1}{a+1}$, $V'(\beta)=\left(\frac{2}{1+a}\right)^{\frac{1}{1-a}}\left(\frac{1-a}{2}\right)$ and $V''(0)=1-a$, we combine this limit with~\eqref{T2} to obtain the equation $\frac{2V(\beta)}{V'(\beta)^2}=\frac{4}{V''(0)}$, equivalent to
\[
(1+a)^{1+a}= 4.
\]
We point out that the previous equality has no solution for $a\in(-1,0)$ since the function $y=x^x$ remains below the constant $1$ if $x\in(0,1)$.
Finally, we consider the case $\alpha\leq -1$. Again $\lim_{x\rightarrow\alpha^+}\frac{\sqrt{2V(x)}}{\abs{V'(x)}}=0$ and  $\lim_{x\rightarrow+\infty}\frac{\sqrt{2V(x)}}{V'(x)}=\frac{T}{\pi}$. On the other hand it is easily checked that $\lim_{x\rightarrow+\infty}\frac{\sqrt{2V(x)}}{V'(x)}=1$ and with a similar reasoning as in the previous case we conclude that $a=-3$.\Qed

There is an alternative proof of the previous result once we know that the potential family must write $V(x)=x-x^{\alpha}$. This second way consists on the computation of the first few Birkhoff coefficients (period constants) on the asymptotic development of the period function at the center. This gives several candidates of parameters $\alpha$ that can be checked to be isochronous potentials by using the characterization with involutions given in~\cite{CMV}. The details of this proof can be found in~\cite[Theorem 3.3]{MRV}. This type of argument is also related to the proofs in~\cite{Albouy,Fejoz}.

\begin{rem}
In the definition of $T$-isochronous family we imposed the uniqueness of the equilibrium. This is essential for the previous proof. To show this we are going to construct a function $\Phi\in C^{\infty}(0,\infty)$ such that the family $V_{\lambda}(x)=\frac{x^2}{2}+\lambda\Phi(x)$ is not $T$-isochronous for any $T$ but every $V_{\lambda}$ is $2\pi$-isochronous around an equilibrium. Let us define $\varphi=\Phi'$ as a function in $C^{\infty}(0,\infty)$ satisfying $\varphi(x)=-3^n$ if $x\in I_n$, $n\in\Z$, where $I_n$ is an open interval containing $[6^n,2\cdot 6^n]$. Note that this function can be constructed because the compact intervals $[6^n,2\cdot 6^n]$ are pairwise disjoint, $2\cdot 6^n<6^{n+1}$. For each $\lambda\in[2^n,2^{n+1}]$ we observe that $\ddot x + x + \lambda\varphi(x)=0$ becomes, on the interval $I_n$,
\[
\ddot x + x - \lambda 3^n=0.
\]
Then $x(\lambda)=\lambda 3^n\in [6^n,2\cdot 6^n]$ is an isochronous center with period $2\pi$.
\end{rem}

\section{Proof of Theorem~\ref{bertrand}}\label{sec:demo}

Let us first be precise on the statement of the Theorem. The attractivity of the force somewhere means that $\phi(\eta)>0$ for some $\eta\in(0,\infty)$. A circular motion is a non-constant solution $(X(t),Y(t))$ of the system~\eqref{central} such that $X(t)^2+Y(t)^2\equiv r_0^2$ is constant. It is easily checked that circular motions have constant angular velocity. Actually, $X(t)+iY(t)=r_0 e^{i\omega t}$ where $\omega=\pm\sqrt{\phi(r_0^2)}$ and $\phi(r_0^2)>0$. The phase space associated to~\eqref{central} can be identified to $(\C\setminus\{0\})\times\C$. In this space we consider the orbit associated to the circular motion
\[
\gamma=\{(r_0\xi,i\omega r_0\xi):\xi\in\mathbb S^1\}.
\]
The assumption on solutions near circular motions means that there exists a neighborhood $\mathcal U$ of $\gamma$, $\mathcal U\subset (\C\setminus\{0\})\times\C$, such that all solutions $(x(t),y(t))$ of \eqref{central} with initial conditions $(x(0)+iy(0),\dot x(0)+i\dot y(0))$ lying in $\mathcal U$ are periodic. In general the size of $\mathcal U$ will depend upon $r_0$.

Identifying the plane of motion with $\C$, setting $\vec{r}=re^{i\theta}$, by Clairaut's change of variable $\rho=r^{-1}$ the central force problem can be written as
\begin{equation}\label{potential}
\frac{d^2 \rho}{d\theta^2} + \rho - \frac{1}{C^2}\rho^{-3}\phi(\rho^{-2})=0,
\end{equation}
where $C=r^2\dot\theta\neq 0$ denotes the angular momentum, which remains constant. Circular motions of~\eqref{central} correspond to equilibrium points of~\eqref{potential}. This type of motion must exist for some values of $C$, indeed we can adjust this parameter so that~\eqref{potential} has an equilibrium on the region where $\phi(\rho^{-2})$ is positive. A well-known fact is the following: A solution of the central potential problem associated to a periodic solution $(\rho(\theta),\rho'(\theta))$ of system~\eqref{potential} is periodic if and only if the minimal period of the solution $(\rho(\theta),\rho'(\theta))$ is commensurable with $\pi$. This period is called the angular period $\Theta$. For future use we present the previous discussion as an auxiliary result.

\begin{lema}\label{lema}
Let $(x(t),y(t))$ be a non-circular periodic solution of~\eqref{central} with $C\neq 0$. Then $\rho=\rho(\theta)$ is a non-constant periodic function whose minimal period is commensurable with $\pi$.
\end{lema}

\begin{demo}
Given a period $T>0$ of the solution of~\eqref{central} we obseve that $r(t+T)=r(t)$ and $\theta(t+T)=\theta(t)+2n\pi$ for some integer $n$. The formula for angular momentum, $C=r^2\dot\theta$ implies that $\dot\theta(t)$ is a $T$-periodic function. Hence the primitive $\theta(t)$ satisfies $\theta(t+T)=\theta(t)+S$ with $S=\int_0^T\frac{C}{r(\tau)^2}d\tau >0$. Hence $S=2n\pi$. Finally we observe that the inverse function $t=t(\theta)$ will satisfy $t(\theta+S)=t(\theta)+T$ and so $\rho(\theta)=\frac{1}{r(t(\theta))}$ has period $S$. The minimal period will be a divisor of $S$; that is, $\frac{2n\pi}{m}$.
\end{demo}

In view of this connection between the system~\eqref{central} and the equation~\eqref{potential} we will employ the following strategy. It will be assumed that the system~\eqref{central} is in the assumptions of Bertrand's Theorem to prove that the family of potentials associated to~\eqref{potential} is $T$-isochronous for some $T>0$ commensurable with $\pi$. Here $\lambda=\frac{1}{C^2}$ and $\Phi'(\rho)=-\rho^{-3}\phi(\rho^{-2})$. This program will be done in three steps.

\begin{claim}
For each $C>0$ the equation~\eqref{potential} has at most one equilibrium.
\end{claim}

We start the proof of this claim with another auxiliary result.

\begin{lema}\label{lema_unicity}
Assume that $V\in C^2(0,\infty)$ has more than one equilibrium. Then one of the following alternatives holds:
\begin{enumerate}[$(i)$]
\item There exists $x_*\in(0,\infty)$ with $V'(x_*)=0$ and a non-constant solution $x(t)$ of $\ddot x + V'(x)=0$ such that $x(t)\rightarrow x_*$ and $\dot x(t)\rightarrow 0$ as $t\rightarrow+\infty$.
\item $E\!:=\{x\in(0,\infty):V'(x)=0\}$ is a non-degenerate interval.
\end{enumerate}
\end{lema}

\begin{demo}
Assume that $(ii)$ does not hold. The set $(0,\infty)\setminus E$ is open and can be partitioned as a disjoint union of open intervals, say $(0,\infty)\setminus E=\cup_{\alpha\in A} I_{\alpha}$ with $A\subset \N$. Assume that some $\alpha$ is such that $I_{\alpha}=(a,b)$ with $0<a<b<\infty$ and $V'(a)=V'(b)=0$, $V'>0$ on $(a,b)$ (alternatively, $V'<0$ on $(a,b)$). We solve $\dot x = +\sqrt{2(V(b)-V(x))}$, $x(0)=\frac{a+b}{2}$ and observe that $x(t)\rightarrow b$ as $t\rightarrow+\infty$. When $V'<0$ on $(a,b)$ we replace the equation by $\dot x = -\sqrt{2(V(a)-V(x))}$. In both cases $(i)$ holds. The property $\dot x(t)\rightarrow 0$ as $t\rightarrow +\infty$ follows from the first order equation.
\end{demo}

After this result is applied to~\eqref{potential}, we are ready to prove the claim. By an indirect argument it is assumed that there are at least two equilibria for some $C>0$. We must look for a contradiction when each of the alternatives given by the Lemma holds.

Assume first that $(i)$ holds. Then there exists a non-constant solution of~\eqref{potential} such that $\rho(\theta)\rightarrow\rho_*$ and $\rho'(\theta)\rightarrow 0$ as $\theta\rightarrow+\infty$. Here $\rho_*>0$ is an equilibrium. The conservation of angular momentum can be interpreted as a first order equation, namely
\[
\dot \theta(t)=C\rho(\theta(t))^2.
\]
Define $m\!:=\inf\{\rho(\theta):\theta\geq\theta(0)\}$ and $M\!:=\sup\{\rho(\theta):\theta\geq\theta(0)\}$. From the assumptions we know that $0<m\leq \dot\theta(t)\leq M<+\infty$ whenever $t\geq 0$ and $\theta(t)$ is well defined. This implies that $\theta(t)$ is well defined on $[0,+\infty)$ and $\theta(t)\rightarrow+\infty$ as $t\rightarrow+\infty$. In consequence the associated solution of~\eqref{central}, $\vec{r}(t)=\frac{1}{\rho(\theta(t))}e^{i\theta(t)}$ is well defined in $[0,\infty)$. Then, as $t\rightarrow+\infty$, $(\vec{r}(t),\dot{\vec{r}}(t))$ accumulates around
\[
\gamma=\{(\tfrac{1}{\rho_*}\xi,C\rho_* i\xi):\xi\in\mathbb S^1\},
\]
an orbit associated to the circular motion $\vec{R}(t)=\frac{1}{\rho_*}e^{iC\rho_*^2 t}$. This is not compatible with the assumption because the solution $\vec{r}(t)$ is not periodic.

To discuss the second alternative we first present another auxiliary result.

\begin{lema}\label{degenera}
Assume that $V\in C^2(0,+\infty)$ and 
\[
E=\{x\in(0,\infty):V'(x)=0\}
\]
is a non-degenerate interval. Let $\gamma_n$ be a sequence of closed orbits of $\dot x=y$, $\dot y=-V'(x)$ approaching $E\times\{0\}$; that is,
\[
\mathrm{dist}(\gamma_n,E\times\{0\})\rightarrow 0 \text{ as }n\rightarrow\infty.
\]
Let $\tau_n>0$ be a sequence of periods associated to $\gamma_n$, then $\tau_n\rightarrow+\infty$.
\end{lema}

In the above statement we have employed the notation
\[
\mathrm{dist}(A,B)=\inf\{\|p-q\| : p\in A, q\in B\}
\]
where $A,B\subset\R^2$.

\begin{demo}
In the plane each closed orbit must surround an equilibrium. Then $\gamma_n$ must go around $E\times\{0\}$, meaning that $E\times\{0\}$ lies in the bounded connected component of $\R^2\setminus\gamma_n$. Assume by contradiction that, after extracting a subsequence $\{\gamma_n\}$, the corresponding periods remain bounded, say $\tau_n\leq \tau <+\infty$. By a compactness argument we can extract a new subsequence $\{\gamma_k\}$ such that $\mathrm{dist}((x_0,0),\gamma_k)\rightarrow 0$, where $x_0$ is some number in $E$. The property of continuous dependence with respect to initial conditions, applied on the interval $[0,\tau]$, implies that $\gamma_k$ converges to $(x_0,0)$ in the Hausdorff distance. This is not possible if $\gamma_k$ goes around $E\times\{0\}$ and the interval $E$ is not a point.
\end{demo}

After this result is applied to the potential $V(\rho)=\frac{1}{2}\rho^2+\frac{1}{C^2}\Phi(\rho)$ we know that the interval $E$ should be of the type $E=[\gamma,\Gamma]$ with $0<\gamma<\Gamma<+\infty$. If $\gamma=0$ or $\Gamma=+\infty$ the orbits of~\eqref{potential} could not go around $E\times\{0\}$ and they could not be closed. This would imply that some non-circular motions of~\eqref{central} could not be periodic and arbitrarily close to circular motions. We will also assume that the condition below holds,
\begin{equation}\label{potmas}
V(\rho)>V(E) \text{ if }\rho\in(0,\infty)\setminus E.
\end{equation}
Otherwise the inequality $V(\rho)<V(E)$ could be valid on $(0,\gamma)$ or $(\Gamma,+\infty)$. The argument of case $(i)$ could be adapted to find non-periodic solutions of~\eqref{central} approaching a circular motion with radius $\tfrac{1}{\gamma}$ or $\tfrac{1}{\Gamma}$. From now on we assume that~\eqref{potmas} holds. The solution of~\eqref{potential} with initial conditions $\rho(0)=\rho_0$, $\rho'(0)=0$ is periodic if $\rho_0\in(\Gamma,\Gamma+\delta)$ and $\delta>0$ is small enough. This solution is denoted by $\rho(\theta,\rho_0)$, with minimal period $\Theta=\Theta(\rho_0)$. We claim that the function $\Theta:(\Gamma,\Gamma+\delta)\rightarrow\R$ is continuous. This can be proved using the formula
\[
\Theta(\rho_0)=\sqrt{2}\int_{A(\rho_0)}^{\rho_0}\frac{d\xi}{\sqrt{V(\rho_0)-V(\xi)}},
\]
where $A$ is the only root of $V(\rho_0)=V(A)$ lying in $(0,\gamma)$. An alternative method, proving that $\Theta$ is indeed in the class $C^1$, could consist in the application of the implicit function Theorem to the equation $\rho'(\Theta,\rho_0)=0$. Going back to Lemma~\ref{degenera} we observe that $\Theta(\rho_0)\rightarrow+\infty$ as $\rho_0\rightarrow\Gamma^+$. This implies that $\Theta$ is not constant and therefore there are many values of $\rho_0\in(\Gamma,\Gamma+\delta)$ such that $\Theta(\rho_0)$ is not commensurable with $\pi$. In this way we construct solutions of~\eqref{central} which are not periodic (in fact they are quasi-periodic) and remain close to the circular solution with radius $\frac{1}{\Gamma}$. The first claim has been proved.

\begin{claim}\label{claim2}
The set $\mathcal C=\{C\in(0,\infty):\eqref{potential} \text{ has an equilibrium}\}$ is connected.
\end{claim}

Define $\Psi(\rho)=\rho^{-4}\phi(\rho^{-2})$ and $\mathcal P=\{\rho\in(0,\infty):\Psi(\rho)>0\}$. We first prove that $\mathcal P$ is an interval. Otherwise there should exist three numbers $0<\rho_1<\rho_2<\rho_3$ with $\Psi(\rho_1)>0$, $\Psi(\rho_2)\leq 0$, $\Psi(\rho_3)>0$. Since $\Psi$ is continuous we can find $\sigma_1\in(\rho_1,\rho_2)$ and $\sigma_2\in(\rho_2,\rho_3)$ with $\Psi(\sigma_1)=\Psi(\sigma_2)>0$. These two numbers would produce two equilibria of~\eqref{potential} for the same $C$. This is against the previous claim.

The set $\mathcal C$ can be expressed as
\[
\mathcal C=\{\sqrt{\Psi(\rho)}:\rho\in\mathcal P\}
\]
and therefore it is connected, proving the second claim.

\begin{claim}\label{claim3}
The family of potentials associated to~\eqref{potential} is $T$-isochronous for some $T\in\pi\Q$.
\end{claim}

At this point it is convenient to emphasize the dependence of the potential associated to~\eqref{potential} with respect to the angular momentum. We write
\[
V(\rho,C)=\frac{1}{2}\rho^2+\frac{1}{C^2}\Phi(\rho).
\]
For each $\rho\in\mathcal C$ we know that $V(\cdot,C)$ has a unique equilibrium producing a circular motion of~\eqref{central}. By the main assumption in Bertrand's Theorem and Lemma~\ref{lema} the equation~\eqref{potential} has a center at $\rho_*=\rho_*(C)$, the critical point of $V(\cdot,C)$. Let $\Theta=\Theta(\rho_0,C)$ be the minimal period of the solution with initial conditions $\rho(0)=\rho_0$, $\rho'(0)=0$. The function $\Theta(\cdot,C)$ is continuous in a neighborhood of $\rho_*(C)$ and $\Theta(\rho_0,C)\in\pi\Q$. This implies that $\Theta(\cdot,C)$ is constant. In consequence $V(\cdot,C)$ is $T(C)$-isochronous around $\rho_*(C)$ with $T(C)=\tfrac{2\pi}{\sqrt{V''(\rho_*(C),C)}}$. We will prove that $T(C)$ is constant. Since $V''(\rho_*(C),C)>0$, the function $\rho_*(C)$ is defined implicitly by the equation $V'(\rho_*(C),C)=0$ and it is $C^1$ on the interval $\mathcal C$. In consequence $T(C)$ is a continuous function defined on $\mathcal C$ and taking values on $\pi\Q$. From Claim~\ref{claim2} we know that $\mathcal C$ is connected and so $T$ is independent of $C$. The family $\{V(\cdot,C)\}_{C>0}$ is $T$-isochronous. This proves the last claim.

On account of Theorem~\ref{fam} the only two possibilities when equilibria exists are $\varphi(\rho)=K$ and $\varphi(\rho)=K\rho^{-3}$. Using the identities $\varphi(\rho)=\rho^{-3}\phi(\rho^{-2})$ and $\rho=r^{-1}$, the corresponding central problems are the Harmonic oscillator $\phi(x)=K$ and Kepler problem $\phi(x)=Kx^{-3/2}$.


\begin{thebibliography}{9}

\bibitem{Albouy}
A.~Albouy,
\emph{Lectures on the two-body problem},
Classical and Celestial Mechanics (Recife, 1993/1999). Princeton University Press, Princeton (2002) 63--116.

\bibitem{Bertrand}
J.~Bertrand,
\emph{Th\'eor\`eme relatif au mouvement d'un point attir\'e vers un centre fixe},
Comptes Rendus Acad. Sci. {\bf 77} (1873) 849--853.

\bibitem{CMV}
A.~Cima, F.~Ma\~nosas, J.~Villadelprat,
\emph{Isochronicity for several classes of Hamiltonian systems},
J. Differential Equations {\bf 157} (1999) 373--413.

\bibitem{Corduneanu}
C.~Corduneanu,
\emph{Principles of differential and integral equations},
AMS Chelsea Pub (1977).

\bibitem{Er1}
V.P.~Ermakov,
\emph{Second order differential equations. Conditions of complete integrability}, Universita Izvestia Kiev, Series III {\bf 9} (1880) 1--25. (Translation from Russian in: Appl. Anal. Discrete Math. {\bf 2} (2008), 123--145.)

\bibitem{Fejoz}
J.~Fejoz, L.~Kaczmarek,
\emph{Sur le th\'eor\`eme de Bertrand (d'apr\`es Michael Herman)},
Ergod. Th. \& Dynam. Sys. {\bf 24} (2004) 1583--1589.

\bibitem{Er2}
P.G.L.~Leach, K.~Andriopoulos,
\emph{The Ermakov Equation: a commentary},  Applicable Analysis and Discrete Mathematics {\bf 2} (2008) 146--157.

\bibitem{MRV}
F.~Ma\~nosas, D.~Rojas, J.~Villadelprat,
\emph{Study of the period function of a two-parameter family of centers},
J. Math. Anal. Appl. {\bf 452} (2017) 188--208.

\bibitem{Pinney}
E.~Pinney,
\emph{The nonlinear differential equation $y''+p(x)y+cy^{-3}=0$},
Proc. Amer. Math. Soc. {\bf 1} (1950) 681.

\bibitem{Urabe61}
M.~Urabe,
\emph{Potential forces which yield periodic motions of a fixed period},
J. Math. Mech. {\bf 10} (1961) 569--578.

\bibitem{Urabe62}
M.~Urabe,
\emph{The potential force yielding a periodic motion whose period is an arbitrary continuous function of the amplitude of the velocity},
Arch. Rational Mech. Anal. {\bf 11} (1962) 27--33.

\bibitem{Rusos}
O.A.~Zagryadskii, E.A.~Kudryavtseva, D.A.~Fedoseev,
\emph{A generalization of Bertrand's theorem to surfaces of revolution},
Sbornik: Mathematics {\bf 203} (2012) 1112--1150.

\end{thebibliography}
\end{document}